\documentclass[letterpaper, 10 pt, conference]{ieeeconf}  % Comment this line out if you need a4paper

\IEEEoverridecommandlockouts                              % This command is only needed if 
\usepackage{cite}
\usepackage{amsmath,amssymb,amsfonts}
\usepackage{algorithmic}
\usepackage{graphicx}
\usepackage{textcomp}
\usepackage{xcolor}
\usepackage{bm}
\usepackage{comment}
\usepackage{lineno,hyperref,float,amsfonts,amsmath,epsfig,graphics,booktabs,multirow,cite,stfloats}

\usepackage{amsthm,mathrsfs}

\theoremstyle{definition}

\theoremstyle{plain}
\newtheorem{myTheo}{Theorem}

\theoremstyle{remark}

\newtheorem{myProp}{Proposition}

\DeclareMathOperator{\diag}{diag}
\DeclareMathOperator{\tr}{tr}
\DeclareMathOperator{\const}{const}

\def\BibTeX{{\rm B\kern-.05em{\sc i\kern-.025em b}\kern-.08em
    T\kern-.1667em\lower.7ex\hbox{E}\kern-.125emX}}
\title{\LARGE \bf Linear quadratic control for discrete-time systems \\with stochastic and bounded noises

}
\author{Xuehui Ma$^{1}$, Shiliang Zhang$^{2}$, Xiaohui Zhang$^{1}$, Jing Xin$^{1}$, Hector Garcia de Marina$^{3}$
\thanks{This work was supported by the National Natural Science Foundation of China (Grant no. 62473311) and by the ERC Starting Grant iSwarm 101076091 and the RYC2020-030090-I grant from the Spanish Ministry of Science. }
\thanks{$^{1}$Xuehui Ma, Xiaohui Zhang and Jing Xin are with the School of Automation and Information Engineering, Xi'an University of Technology, Xi'an, China ({\tt\small xuehui.yx@gmail.com, xhzhang@xaut.edu.cn, xinj@xaut.edu.cn}).}
\thanks{$^{2}$Shiliang Zhang is with the Department of Informatics, University of Oslo, Norway ({\tt\small shilianz@ifi.uio.no}).} 
\thanks{$^{3}$Hector Garcia de Marina is with the Department of Computer Engineer, Automation and Robotics, and with CITIC, University of Granada, Spain ({\tt\small hgdemarina@ugr.es}).}
}

\begin{document}
\maketitle

\begin{abstract}
This paper focuses on the linear quadratic control (LQC) design of systems corrupted by both stochastic noise and bounded noise simultaneously. When only of these noises are considered, the LQC strategy leads to stochastic or robust controllers, respectively. However, there is no LQC strategy that can simultaneously handle stochastic and bounded noises efficiently. This limits the scope where existing LQC strategies can be applied. In this work, we look into the LQC problem for discrete-time systems that have both stochastic and bounded noises in its dynamics. We develop a state estimation for such systems by efficiently combining a Kalman filter and an ellipsoid set-membership filter. The developed state estimation can recover the estimation optimality when the system is subject to both kinds of noise, the stochastic and the bounded. Upon the estimated state, we derive a robust state-feedback optimal control law for the LQC problem. The control law derivation takes into account both stochastic and bounded-state estimation errors, so as to avoid over-conservativeness while sustaining stability in the control. In this way, the developed LQC strategy extends the range of scenarios where LQC can be applied, especially those of real-world control systems with diverse sensing which are subject to different kinds of noise. We present numerical simulations, and the results demonstrate the enhanced control performance with the proposed strategy.
\end{abstract}

\section{Introduction} \label{section1}

%%%% why stochastic and bounded noises?  background? practical applications?
In this work, we study the linear quadratic control (LQC) problem for discrete-time systems under process and observation noise. The construction of a control strategy for this LQC problem depends on how the noise is modeled. In general, noise is modeled as (i) stochastic described by its probability distribution \cite{11219403,ma2022adaptivequantile} and (ii) bounded, depicted by compact sets that confine all possible values for the noise \cite{10591250,lorenzen2019robust}. In the case of stochastic noise, the classical linear-quadratic-Gaussian (LQG) control strategies can be applied to optimize the control~\cite{aastrom2012introduction}. On the other hand, for bounded noise, min-max control approaches can provide the optimal robust control law~\cite{bertsimas2007constrained}.

However, the two categories of noise modeling and the construction of a control strategy typically work in parallel, and each of the two categories assumes only one type of noise in its control design. However, this assumption does not hold in real world scenarios~\cite{di2004set,hanebeck1999new,liu2025adaptive}. Taking autonomous vehicles as an example, both stochastic and bounded noises are present in vehicle path planning and obstacle avoidance~\cite{SOLIMAN2022387,yang2022set,ma2025fault,liu2025adaptiveFault}. \textit{I.e.}, the noise induced by inertia sensors, LiDAR and camera in a vehicle is stochastic noise that can be described by statistical models, while the noise in the tire grip and actuators can be better modeled by a bounded set. The existence of both stochastic and bounded noises in the same system requires the design of an LQC strategy that can deal with both kinds of noise simultaneously. 

%%%%% the limitations of statistic thoery based estimator and controller (e.g., Kalman filter and LQG controller) (estimation bias, over-optimistic, no robustness)
The existing two categories of LQC strategies - either stochastic control or robust control - cannot work efficiently when there are stochastic and bounded noises. Despite the celebrated success of stochastic control for systems with stochastic noise, such approaches are limited in their ability to handle bounded noise~\cite{MISHRA2021109512,arcari2023stochastic,ma2022active}. Most of the stochastic control approaches use the Kalman filter to estimate system states, with the assumption that the process and observation noise are Gaussian distributed. Unfortunately, such an assumption no longer holds in the case of bounded noise, leading to an overoptimistic estimation that goes into the controller ~\cite{de2024state,lu2017multi}. Such an overoptimism produces increased state estimation errors, and the propagation of the errors to the control loop degrades the control performance. As a result, there is a risk of divergence from the actual values and their desired trajectories, since the estimation errors accumulate over time.

%%%%% the limitations of set thoery based estimator and controller (e.g., set-membership filter and min-max controller) (over-conservativeness )
Unlike stochastic control that relies on the statistic distribution of noise, robust control only needs the boundary information of the noise in its control law derivation~\cite{bertsimas2007constrained}. Robust control approaches, like min-max control, seek robustness against the most pessimistic noise within the bounded set that has a worst impact on the control. However, the performance of robust control becomes over-conservative in the presence of stochastic noise~\cite{10591250,lorenzen2019robust}. Robust control guarantees robustness against all possible noise, including the extreme noise that rarely occurs. The consideration of such extreme noise reduces the control performance when the system works under normal conditions, leading to very conservative controllers. Recent studies integrated set-membership filters into robust control to learn the noise boundary and reduce its conservativeness~\cite{iannelli2020structured,parsi2022explicit,parsi2022scalable,10551447}. However, these studies neglect the case where the stochastic noise exceeds the learned boundary, making the set-membership filter infeasible for the controller~\cite{10904011,noack2012combined}. 

%%%% this work ..... and contributions
This paper aims to address the challenges in linear quadratic control for systems with the presence of both stochastic and bounded noises. To this end, we develop the state estimation and derive the LQC control law, where the stochastic and bounded noises are incorporated simultaneously. Specifically, we propose a mixed-noise model that unifies the mathematical description of stochastic and bounded noises. We then derive the state estimation based on the mixed-noise model. The derived state estimation considers both kinds of noise, and leads to an estimation that is less optimistic than the produced by the Kalman filter, while it is still less conservative than the set-membership filter. Finally, we integrate the state estimation in the derivation of a robust control law, where both the stochastic and bounded estimation errors are considered and minimized. Simultaneously minimizing these two estimation errors avoids over-conservativeness and sustains stability in the control. We demonstrate how our approach improves the estimation and control performance compared with robust control with Kalman filter and set-membership filter. 

This paper is organized as follows.	We formulate the linear quadratic control problem with stochastic and bounded noises in Section~\ref{section2}. Section \ref{section3} derives the state estimation that incorporates the Kalman filter and the set-membership filter, and Section~\ref{section4} presents the optimal control law derivation. In Section~\ref{section5}, we present numerical simulations to demonstrate the effectiveness of our approach. Finally, we conclude our work in Section~\ref{section6}.

\textit{Notations:} 

We denote the ellipsoid set $E(\bm{c},\bm{M})$ in the set-membership approach as, 
\begin{equation}\label{defin_ellip}
    \mathcal{E}(\bm{c},\bm{M})= \{ \bm{x} \in \mathbb{R}^n: (\bm{x}-\bm{c})^T \bm{M}^{-1} (\bm{x}-\bm{c}) \leq 1 \},
\end{equation}
where $\bm{c}$ is the center of the ellipsoid and $\bm{M}$ is the \emph{shape matrix} of the ellipsoid. The \emph{Minkowski sum} of two sets of vectors $\bm{A}$ and $\bm{B}$ in Euclidean space is defined by adding each vector in $\bm{A}$ to each vector in $\bm{B}$, written as
\begin{equation}
\bm{A} \oplus \bm{B} = \{ \bm{a}+\bm{b} \, | \, \bm{a} \in \bm{A}, \bm{b}\in \bm{B} \}.
\end{equation}

\section{Problem Formulation}\label{section2}
In this work, we consider the following discrete-time linear dynamic system 
\begin{equation}\label{systemdynamic}  \bm{x}_{k+1}=\bm{A}_{k}\bm{x}_{k}+\bm{B}_{k}\bm{u}_{k}+\bm{w}_{k}, 
\end{equation}
\begin{equation}\label{observationequation}
    \bm{z}_{k} = \bm{H}_{k}\bm{x}_{k}+\bm{v}_{k}, 
\end{equation}
where $\bm{x}_{k} \in \mathbb{R}^n$ is the state vector and cannot be measured directly, $\bm{z}_{k} \in \mathbb{R}^m$ is the measurement, $\bm{u}_{k} \in \mathbb{R}^r$ is the control vector, $\bm{A}_k \in \mathbb{R}^{n\times n}$ is the state matrix, $\bm{B}_k\in \mathbb{R}^{n\times r}$ is the control matrix, and $\bm{H}_k\in \mathbb{R}^{m\times n}$ is the measurement matrix. $\bm{w}_{k} \in \mathbb{R}^n$ and $\bm{v}_k \in \mathbb{R}^m$ are process and measurement noises. Note that in real-world applications, process noise $\bm{w}_{k}$ and measurement noise $\bm{v}_k$ typically consist of both stochastic and bounded noises. For example, in the operation of autonomous wheeled mobile robot, the process noise $\bm{w}_{k}$ induced by motor drive current fluctuations is often modeled as Gaussian, whereas external collisions or pushing disturbances are better represented by bounded sets. Similarly, the measurement noise $\bm{v}_k$ in the mobile robot operation may include Gaussian components like IMU noise and GPS jitter, and also bounded components, such as IMU mis-calibration and GPS position error caused by multi-path effects.

Therefore, we incorporate both stochastic and bounded noises to better represent system uncertainties. Particularly, we model the process and observation noises in \eqref{systemdynamic} and \eqref{observationequation} by superposing both stochastic and bounded noises
\begin{equation}\label{noise_w}
    \bm{w}_k = \bm{w}_k^s +\bm{w}_k^b,
\end{equation}
\begin{equation}\label{noise_v} 
    \bm{v}_k = \bm{v}_k^s +\bm{v}_k^b,
\end{equation}
where $\bm{w}_k^s$ and $\bm{v}_k^s$ are stochastic noises, $\bm{w}_k^b$ and $\bm{v}_k^b$ are bounded noises. Here we assume the stochastic noises follow  Gaussian distributions:
\begin{equation}\label{stoc_w}
    \bm{w}_k^s \sim \mathcal{N}(0,\bm{P}_k^w),  
\end{equation}
\begin{equation}\label{stoc_v}
    \bm{v}_k^s \sim \mathcal{N}(0,\bm{P}_k^v), 
\end{equation}
where the noises are white with their mean equal to zero, and $\bm{Q}_k^w$ and $\bm{Q}_k^v$ are the associated covariances. 
The bounded noises are assumed to be bounded within the following ellipsoidal sets, respectively:
\begin{equation}\label{bound_w}
    \bm{w}_k^b \in \mathcal{E}(0,\bm{M}_k^w),
\end{equation}
\begin{equation}\label{bound_v} 
    \bm{v}_k^b \in \mathcal{E}(0,\bm{M}_k^v),
\end{equation}
where $\bm{M}_k^w$ and $\bm{M}_k^v$ are the matrices defining the shape, size, and orientation of the ellipsoids. The ellipsoidal centers are zero. We model the initial state vector $\bm{x}_0$ for the noise as 
\begin{equation}\label{state_x0}
\begin{aligned}
    &\bm{x}_0 = \bar{\bm{x}}_0 +\bm{e}_0^s + \bm{e}_0^b, \\
\end{aligned}
\end{equation}
where $\bar{\bm{x}}_0$ is the center of the system state, $\bm{e}_0^s$ and $\bm{e}_0^b$ are stochastic and bounded uncertainties due to stochastic and bounded noises, respectively. The distribution of $\bm{e}_0^s$ and the boundary of $\bm{e}_0^b$ are:
\begin{equation}\label{state_x0_prior}
\begin{aligned}
    &\bm{e}_0^s \sim \mathcal{N}(0,\bm{P}_0),\\
    &\bm{e}_0^b \in \mathcal{E}(0,\bm{M}_0),
\end{aligned}
\end{equation}
where $\bm{P}_0$ is the covariance of Gaussian distribution, and $\bm{M}_0$ is the shape matrix of the ellipsoid set.

With the settings above, this work aims to control the system \eqref{systemdynamic}-\eqref{observationequation} that seeks minimizing the following quadratic cost function
\begin{equation} \label{J_costF}
\begin{aligned}
    J = & \bm{x}^T_N\bm{Q}_N\bm{x}_N +\sum_{k= 0}^{N-1} \left( \bm{x}^T_k\bm{Q}_k\bm{x}_k +\bm{u}^T_k\bm{R}_k\bm{u}_k\right) , 
\end{aligned}
\end{equation}
where the state-cost matrix $\bm{Q}_k \in \mathbb{R}^{n\times n}$ and control-cost matrix $\bm{R}_k \in \mathbb{R}^{r\times r}$ are positive semidefinite and positive definite, respectively. 

Our hypothesis is that the process noise, observation noise, and initial state can be jointly described by a stochastic model and a bounded set. In this way, the control will be different from LQG control or min-max control, \textit{i.e.}, instead of minimizing the expectation of the cost function as in LQG or the worst value of cost function as in min-max control, we minimize the expected value and the worst value for the stochastic noise and bounded set noises in the cost function $J$ \eqref{J_costF}, respectively. 
We formulate this control optimization problem as
\begin{equation}\label{Prob1}
\begin{aligned}
    (\mathcal{P}): \quad  &\min_{\left\{\bm{u}_k\vert_{k=0}^{N-1}\right\}} \quad \max_{\{\bm{w}_k^b,\bm{v}_k^b,\bm{e}_0^b\}} \left\{\mathbb{E}_{\{\bm{w}_k^s, \bm{v}_k^s, \bm{e}_0^s\}} \left\{ J \right\} \right\}\\
    \quad s.t. \quad &\bm{x}_{k+1}=\bm{A}_{k}\bm{x}_{k}+\bm{B}_{k}\bm{u}_{k}+\bm{w}_{k},\\
    \quad \quad \quad \quad & \bm{z}_{k} = \bm{H}_{k}\bm{x}_{k}+\bm{v}_{k},  k = 0,1,\cdots,N-1\\
    \quad \quad \quad \quad & \bm{w}_k^s \sim \mathcal{N}(0,\bm{P}_k^w),\\
    \quad \quad \quad \quad & \bm{w}_k^b \in \mathcal{E}(0,\bm{M}_k^w),\\
    \quad \quad \quad \quad & \bm{v}_k^s \sim \mathcal{N}(0,\bm{P}_k^v),\\ 
    \quad \quad \quad \quad & \bm{v}_k^b \in \mathcal{E}(0,\bm{M}_k^v),\\
    \quad \quad \quad \quad & \bm{e}_0^s \sim \mathcal{N}(0,\bm{P}_0),\\
    \quad \quad \quad \quad & \bm{e}_0^b \in \mathcal{E}(0,\bm{M}_0),
\end{aligned}
\end{equation}
and the optimal control sequence $\bm{u}_k\vert_{k=0}^{N-1}$ can be derived by solving the optimization problem $(\mathcal{P})$. 

As defined in \eqref{J_costF}, the cost $J$ in $(\mathcal{P})$ is a function of the system state $\bm{x}_k$. When combining the equation \eqref{systemdynamic}, the cost $J$ can be rewritten as the function of initial state $\bm{x}_0$. Then we can gain the optimal control sequence $\bm{u}_0, \cdots, \bm{u}_{N-1}$ by minimizing $\max_{\{\bm{w}_k^b,\bm{v}_k^b,\bm{e}_0^b\}} \left\{\mathbb{E}_{\{\bm{w}_k^s, \bm{v}_k^s, \bm{e}_0^s\}} \left\{ J(\bm{x}_0) \right\} \right\}$ in problem \eqref{Prob1}. However, the gained control sequence is overconservative, because the boundary of error $\bm{e}_0^b$ for the initial state is usually set excessively wide to cover and tolerate the worst case. This results in the over-conservatism in control, since the worst case rarely occurs. In addition, the covariance of stochastic error $\bm{e}_0^s$ is also initiated with excessively large values, which cannot quantify the stochastic uncertainty accurately. 

In this paper, we will use the measurements $\bm{z}_{k}$ to estimate and update the distribution and boundary of the states, and we present such a state estimation in detail in Section~\ref{section3}. Note that the update of stochastic and bounded uncertainties in system states complicates the optimization problem in $(\mathcal{P})$, due to the high complexity calculation of the expectation and the maximum value in the cost $J$ in $(\mathcal{P})$. To address this issue, we will reformulate problem $(\mathcal{P})$ and make it tractable to derive the optimal control law. We detail the reformulation of problem $(\mathcal{P})$ and control law derivation in Section~\ref{section4}.

\section{State Estimation incorporating Kalman Filter and Set-membership Filter}\label{section3}

This section derives a recursive state estimation for the system \eqref{systemdynamic}-\eqref{observationequation} with the presence of both stochastic and bounded noises. In Section \ref{section2}, we make the assumption that the stochastic noise follows a Gaussian distribution and the bounded noises are within an ellipsoid set. Although Kalman and ellipsoidal set-membership filters are effective for white Gaussian noise and ellipsoid set bounded noises, respectively, these two filters are no longer applicable for estimating the stochastic and bounded uncertainties simultaneously. Therefore, we propose the combination of Kalman and ellipsoidal set-membership filters for the design of the optimal state estimation. We construct the estimation by minimizing the covariance of the stochastic estimation error and the ellipsoid volume of the bounded estimation error simultaneously~\cite{maksarov1996state,noack2012optimal}. 

Let us show the derivation of our recursive state estimation. We start by defining the estimation of state at the $k$-th instant given the measurements $\bm{z}_{k}$ as 
\begin{equation}\label{state_x_hat}
   \bm{x}_k \sim \mathcal{X}(\hat{\bm{x}}_{k|k}, \bm{P}_{k|k}, \bm{M}_{k|k}),
\end{equation}
where $\hat{\bm{x}}_{k|k}$ is the mixed center of the estimation with stochastic and bounded uncertainties, $\bm{P}_{k|k}$ is the covariance matrix of stochastic estimation error, $\bm{M}_{k|k}$ is the shape matrix of ellipsoid set estimation error. Combining the main technique in the Kalman and the ellipsoidal set-membership filters, we update the $\hat{\bm{x}}_{k|k}$, $\bm{P}_{k|k}$ and $\bm{M}_{k|k}$ by utilizing the measurements $\bm{z}_{k}$, as shown in the following formal statement. 

\begin{myProp}\label{Theo1}
Consider the discrete-time linear system described by \eqref{systemdynamic} and \eqref{observationequation}, which is disturbed by stochastic and bounded noises as in (\ref{noise_w})-(\ref{bound_v}). Let the initial state be set as \eqref{state_x0} and \eqref{state_x0_prior}, with initial mixed state center $\bar{\bm{x}}_0$, initial covariance $\bm{P}_0$, and initial shape matrix $\bm{M}_0$. Then we update $\hat{\bm{x}}_{k|k}$, $\bm{P}_{k|k}$ and $\bm{M}_{k|k}$ recursively in the following two steps:\\
(i) Predict
\begin{equation}\label{timeupdate1}
\begin{aligned}
    \hat{\bm{x}}_{k|k-1}=\bm{A}_{k-1}\hat{\bm{x}}_{k-1|k-1}+\bm{B}_{k-1}\bm{u}_{k-1},
\end{aligned}
\end{equation}
\begin{equation}\label{timeupdate2}
\begin{aligned}
    \bm{P}_{k|k-1}=\bm{A}_{k-1}\bm{P}_{k-1|k-1}\bm{A}^T_{k-1}
    +\bm{P}^w_{k-1},
\end{aligned}
\end{equation}
\begin{equation}\label{timeupdate3}
\begin{aligned}
    \bm{M}_{k|k-1}=&(p^{-1}_k+1)\bm{A}_{k-1}\bm{M}_{k-1|k-1}\bm{A}^T_{k-1}\\
    &+(p_k+1)\bm{M}^w_{k-1},
\end{aligned}
\end{equation}
(ii) Update
\begin{equation}\label{obserupdate1}
    \hat{\bm{x}}_{k|k}=\hat{\bm{x}}_{k|k-1}+\bm{\Gamma}_k\left( \bm{z}_k-\bm{H}_k\hat{\bm{x}}_{k|k-1} \right) ,
\end{equation}
\begin{equation}\label{obserupdate2}
    \bm{P}_{k|k}=(\bm{I}-\bm{\Gamma}_k\bm{H}_k)\bm{P}_{k|k-1}(\bm{I}-\bm{\Gamma}_k\bm{H}_k)^T +\bm{\Gamma}_k\bm{P}_k^v\bm{\Gamma}_k^T,
\end{equation}
\begin{equation}\label{obserupdate3}
\begin{aligned}
    \bm{M}_{k|k}=&(q^{-1}_k+1)(\bm{I}-\bm{\Gamma}_k\bm{H}_k)\bm{M}_{k|k-1}(\bm{I}-\bm{\Gamma}_k\bm{H}_k)^T\\
    &+(q_k+1)\bm{\Gamma}_k\bm{M}^v_{k}\bm{\Gamma}_k^T,
\end{aligned}
\end{equation}
\begin{equation}\label{KalmanG}
\begin{aligned}
    \bm{\Gamma}_k = &\left[ \bm{P}_{k|k-1}\bm{H}_k^T + (q^{-1}_k+1)\bm{M}_{k|k-1}\bm{H}_k^T \right] \\
    &\left[ (q^{-1}_k+1)\bm{H}_k\bm{M}_{k|k-1}\bm{H}_k^T +(q_k+1)\bm{M}^v_{k} \right.\\
    &\left.  + \bm{H}_k\bm{P}_{k|k-1}\bm{H}_k^T+\bm{P}_k^v \right]^{-1},
\end{aligned}
\end{equation}
where the scalar parameters $p_k$ and $q_k$ are calculated by 
\begin{equation}\label{p_k}
p_k = \left( \frac{\tr(\bm{A}_{k-1}\bm{M}_{k-1|k-1}\bm{A}^T_{k-1})}{\tr(\bm{M}^w_{k})}\right)^{\frac{1}{2}},
\end{equation}
and 
\begin{equation}\label{q_k}
q_k = \left(\frac{\tr((\bm{I}-\bm{\Gamma}_k)\bm{M}_{k|k-1}(\bm{I}-\bm{\Gamma}_k)^T)}{\tr(\bm{\Gamma}_k\bm{M}^w_{k}\bm{\Gamma}_k^T)}\right)^{\frac{1}{2}}.
\end{equation}
\end{myProp}

\begin{proof}
Combine the system defined in \eqref{systemdynamic} and the process noise defined in \eqref{noise_w}, and we have
\begin{equation}\label{dynamic_x_k}
    \bm{x}_{k}=\bm{A}_{k-1}\bm{x}_{k-1}+\bm{B}_{k-1}\bm{u}_{k-1}+\bm{w}_{k-1}^s +\bm{w}_{k-1}^b,
\end{equation}
where $\bm{w}_{k-1}^s \sim \mathcal{N}(0,\bm{P}_{k-1}^w)$ and $\bm{w}_{k-1}^b \in \mathcal{E}(0,\bm{M}_{k-1}^w)$. We use the mixed center $\hat{\bm{x}}_{k|k}$ defined in \eqref{state_x_hat} to represent the estimated value of the state at the instant $k$. Then we define the estimated error as
\begin{equation}\label{estimate_error_k1}
    \bm{e}_k = \bm{x}_k-\hat{\bm{x}}_{k|k},
\end{equation}
where the estimated error consists of a stochastic part and a bounded part, shown as
\begin{equation}\label{estimate_error_k2}
\begin{aligned}
    &\bm{e}_k = \bm{e}_k^s + \bm{e}_k^b,\\
    &\bm{e}_k^s \sim \mathcal{N}(0,\bm{P}_k),\\
    &\bm{e}_k^b \in \mathcal{E}(0,\bm{M}_k).
\end{aligned}
\end{equation}
We use the system dynamics in \eqref{dynamic_x_k} to predict the mixed center $\hat{\bm{x}}_{k|k-1}$, as shown in \eqref{timeupdate1}. The predicted covariance matrix is calculated by 
\begin{equation}
    \bm{P}_{k|k-1}= \mathbb{E}(\bm{A}_{k-1}\bm{e}_{k-1}^s+\bm{w}_{k-1}^s)^2,
\end{equation}
and then we get \eqref{timeupdate2}. We utilize the Minkowski sum to calculate the predicted shape matrix of the ellipsoid set 
\begin{equation}
\begin{aligned}
    &\bm{A}_{k-1}\bm{e}_{k-1}^b+\bm{w}_{k-1}^b \\
    \in &\bm{A}_{k-1}\mathcal{E}(0,\bm{M}_{k-1|k-1}) \oplus \mathcal{E}(0,\bm{M}^w_{k-1})\\
    \subset & \mathcal{E}(0,\bm{M}_{k|k-1}),
\end{aligned}
\end{equation}
and we obtain \eqref{timeupdate3}. 
Then, we use the measurement data $\bm{z}_{k}$ to update the mixed center \eqref{obserupdate1}, where the gain $\bm{\Gamma}_k$ is obtained by minimizing the updated estimation error $\bm{e}_k$. Substituting the update the mixed center \eqref{obserupdate1} into the error \eqref{estimate_error_k1} yields 
\begin{equation}\label{error_1}
\begin{aligned}
    \bm{e}_k = \bm{x}_k - \hat{\bm{x}}_{k|k-1}- \bm{\Gamma}_k\left( \bm{z}_k-\bm{H}_k\hat{\bm{x}}_{k|k-1} \right), 
\end{aligned}
\end{equation}
With the measurement $\bm{z}_{k} = \bm{H}_{k}\bm{x}_{k}+\bm{v}_k^s +\bm{v}_k^b$ substituted to \eqref{error_1}, we obtain  
\begin{equation}\label{error_2}
\begin{aligned}
    \bm{e}_k = (\bm{I}-\bm{\Gamma}_k\bm{H}_k)(\bm{x}_k - \hat{\bm{x}}_{k|k-1})- \bm{\Gamma}_k\left( \bm{v}_k^s +\bm{v}_k^b \right). 
\end{aligned}
\end{equation}
We define predicted estimation error as $\bm{x}_k - \hat{\bm{x}}_{k|k-1}=\bm{e}_{k|k-1}$. Then, we rewrite the error in \eqref{error_2} as
\begin{equation}\label{error_3}
    \bm{e}_k = (\bm{I}-\bm{\Gamma}_k\bm{H}_k)(\bm{e}_{k|k-1}^s+\bm{e}_{k|k-1}^b)- \bm{\Gamma}_k\left( \bm{v}_k^s +\bm{v}_k^b \right)\\
\end{equation}
According to \eqref{error_3}, we update the covariance matrix by
\begin{equation}
    \bm{P}_{k|k}= \mathbb{E}[(\bm{I}-\bm{\Gamma}_k\bm{H}_k)\bm{e}_{k|k-1}^s -\bm{\Gamma}_k\bm{v}_k^s]^2,
\end{equation}
which leads to \eqref{obserupdate2}. We update the shape matrix of the ellipsoid set as
\begin{equation}
\begin{aligned}
    &(\bm{I}-\bm{\Gamma}_k\bm{H}_k)\bm{e}_{k|k-1}^b -\bm{\Gamma}_k\bm{v}_k^b\\
    \in &(\bm{I}-\bm{\Gamma}_k\bm{H}_k)\mathcal{E}(0,\bm{M}_{k|k-1}) \oplus \bm{\Gamma}_k\mathcal{E}(0,\bm{M}^v_{k})\\
    \subset & \mathcal{E}(0,\bm{M}_{k|k}),
\end{aligned}
\end{equation}
and get \eqref{obserupdate3}.
For the simplicity of the calculation, we derive the optimal $\bm{\Gamma}_k$ by minimizing the trace of covariance matrix $\bm{P}_{k|k}$ and shape matrix $\bm{M}_{k|k}$, and define the estimation cost function as 
\begin{equation}\label{error_4}
\begin{aligned}
    V_k = &\tr\{\bm{P}_{k|k}\}+ \tr\{ \bm{M}_{k|k}\}\\
     =  &\tr\{ (\bm{I}-\bm{\Gamma}_k\bm{H}_k)\bm{P}_{k|k-1}(\bm{I}-\bm{\Gamma}_k\bm{H}_k)^T \} \\
     & +\tr\{\bm{\Gamma}_k\bm{P}_k^v\bm{\Gamma}_k^T \}\\
    & + (q^{-1}_k+1)\tr\{ (\bm{I}-\bm{\Gamma}_k\bm{H}_k)\bm{M}_{k|k-1}(\bm{I}-\bm{\Gamma}_k\bm{H}_k)^T \}\\
    & + (q_k+1)\tr\{\bm{\Gamma}_k\bm{M}^v_{k}\bm{\Gamma}_k^T\}.
\end{aligned}
\end{equation}
Let $\partial V_k /\partial \bm{\Gamma}_k = 0$, we get the optimal gain $\bm{\Gamma}_k$.
\end{proof}

\section{Optimal Controller Design}\label{section4}

In this section, we integrate the iterative estimated states described in Section \ref{section3} into the derivation of an optimal linear quadratic control law. The estimator updates both stochastic and bounded estimation errors of the states at each instant, where the stochastic error $\bm{e}_k^s$ is calculated by the expectation and the bounded error $\bm{e}_k^b$ is added as a constraint in the optimal control problem formulation. 

Let us show how we formulate the optimization problem by defining the quadratic cost function at the $k$-th instant as
\begin{equation}\label{J_k}
    J_k =  \sum_{t=k}^{k+N-1}\bm{x}^T_{t+1}\bm{Q}_{t+1}\bm{x}_{t+1} + \bm{u}^T_{t}\bm{R}_{t}\bm{u}_{t}.
\end{equation}
With the state estimation integrated, the linear quadratic optimal control problem at the $k$-th instant is presented as
\begin{equation}\label{Prob_k}
\begin{aligned}
    (\mathcal{P}_k):  &\min_{\left\{\bm{u}_k\right\}}  \max_{\{\bm{w}_k^b,\bm{v}_k^b,\bm{e}_k^b\}} \left\{ \mathbb{E}_{\{\bm{w}_k^s, \bm{v}_k^s,{e}_k^s\}} \left\{ J_k  \right\} \right\}\\
    \quad s.t. \quad &\bm{x}_{k+1}=\bm{A}_{k}\bm{x}_{k}+\bm{B}_{k}\bm{u}_{k}+\bm{w}_{k},\\
    \quad \quad \quad \quad & \bm{z}_{k} = \bm{H}_{k}\bm{x}_{k}+\bm{v}_{k},  \\
    \quad \quad \quad \quad & \bm{w}_k^s \sim \mathcal{N}(0,\bm{P}_k^w),\\
    \quad \quad \quad \quad & \bm{w}_k^b \in \mathcal{E}(0,\bm{M}_k^w),\\
    \quad \quad \quad \quad & \bm{v}_k^s \sim \mathcal{N}(0,\bm{P}_k^v),\\ 
    \quad \quad \quad \quad & \bm{v}_k^b \in \mathcal{E}(0,\bm{M}_k^v),\\
    \quad \quad \quad \quad & \bm{e}_k^s \sim \mathcal{N}(0,\bm{P}_k),\\
    \quad \quad \quad \quad & \bm{e}_k^b \in \mathcal{E}(0,\bm{M}_k).
\end{aligned}
\end{equation}
We detail the derivation of the control law $\bm{u}_k$ by solving the problem $(\mathcal{P}_k)$ in the following technical result. 

\begin{myProp}\label{Prop1}
The expectation mean of the cost function $J_k$ in $(\mathcal{P}_k)$ can be written in the form
\begin{equation}\label{J_K_exp}
\begin{aligned}
     \mathbb{E}_{\{\bm{w}_k^s, \bm{v}_k^s,{e}_k^s\}} \{J_k\} =&\hat{\bm{x}}_{k}^T\mathscr{A}_k\hat{\bm{x}}_{k} +\bm{U}_{k}^T\mathscr{B}_k\bm{U}_{k} \\
      &+(\bm{\eta}_k^b)^T\mathscr{C}_k\bm{\eta}_k^b+2\hat{\bm{b}}_k^T\bm{U}_{k}  \\
      &+2\bm{U}_{k}^T\mathscr{D}_k\bm{\eta}_k^b+2\hat{\bm{c}}_k^T\bm{\eta}_k^b +\const,\\
\end{aligned}
\end{equation}
for appropriate matrices $\mathscr{A}_k\in\mathbb{R}^{n\times n}$, $\mathscr{B}_k\in\mathbb{R}^{N\cdot r\times N\cdot r}$, $\mathscr{C}_k\in\mathbb{R}^{2N\cdot n\times 2N\cdot n}$, $\mathscr{D}_k\in\mathbb{R}^{N\cdot r\times 2N\cdot n}$, vectors $\hat{\bm{b}}_k\in\mathbb{R}^{N\cdot r}$, $\hat{\bm{c}}_k\in\mathbb{R}^{2N\cdot n}$, and a constant $\const$. The control vector $\bm{U}_k$ and noise vector $\bm{W}_k$ are defined as
\begin{equation}
     \bm{U}_k^T = \left[ \bm{u}_k^T \quad \bm{u}_{k+1}^T \quad \cdots \quad \bm{u}_{k+N-1}^T \right], 
\end{equation}
and
\begin{equation}
     \bm{W}_k^T = \left[ \bm{w}_k^T \quad \bm{w}_{k+1}^T \quad \cdots \quad \bm{w}_{k+N-1}^T \right].
\end{equation}
\end{myProp}

\begin{proof}
With the linear dynamic system \eqref{systemdynamic}, we write the state at the instant $t+1$ as
\begin{equation} \label{state_tplus1}
     \bm{x}_{t+1} = \tilde{\bm{A}}_{t} \bm{x}_k+\tilde{\bm{B}}_{t}\bm{U}_k+\bm{C}_{t}\bm{W}_k,
\end{equation}
where
\begin{equation}
     \tilde{\bm{A}}_{t}=\prod_{i=t}^{k}\bm{A}_i,
\end{equation}
\begin{equation}
     \tilde{\bm{B}}_{t} = \left[ \left( \prod_{i=t}^{k+1}\bm{A}_i \right)\bm{B}_k \quad \cdots \quad \bm{B}_{t} \quad \bm{0}_{n \times (k+N-t-1)\cdot r} \right], 
\end{equation}
\begin{equation}
     \bm{C}_{t} = \left[  \prod_{i=t}^{k+1}\bm{A}_i  \quad \cdots \quad \bm{I}_{n\times n} \quad \bm{0}_{n \times (k+N-t-1)\cdot r} \right], 
\end{equation}
Substituting the estimated errors \eqref{estimate_error_k1} and \eqref{estimate_error_k2} into \eqref{state_tplus1}, we get
\begin{equation}\label{state_form1}
\bm{x}_{t+1} = \tilde{\bm{A}}_{t} (\hat{\bm{x}}_k+\bm{e}_k^s+\bm{e}_k^b) +\tilde{\bm{B}}_{t}\bm{U}_k+\bm{C}_{t}(\bm{W}_k^s+\bm{W}_k^b).
\end{equation}
Define the matrix $\tilde{\bm{C}}_{t}$, stochastic vector $\bm{\eta}_k^s$ and bounded uncertainties vectors $\bm{\eta}_k^b$ as
\begin{equation}\label{Para_C_eta}
\begin{aligned}
    \tilde{\bm{C}}_{t}= \left[ \begin{array}{cc} \tilde{\bm{A}}_{t} & \bm{C}_{t} \end{array} \right],\quad \bm{\eta}_k^s = \left[ \begin{array}{c} \bm{e}_k^s \\ \bm{W}_k^s \end{array} \right], \quad \bm{\eta}_k^b = \left[ \begin{array}{c} \bm{e}_k^b \\ \bm{W}_k^b \end{array} \right],
\end{aligned}
\end{equation}
With \eqref{Para_C_eta}, equation \eqref{state_form1} can be rewritten as
\begin{equation}\label{state_form2}
    \bm{x}_{t+1} = \tilde{\bm{A}}_{t}\hat{\bm{x}}_{k} +\tilde{\bm{B}}_{t}\bm{U}_k +\tilde{\bm{C}}_{t}\bm{\eta}_k^b +\tilde{\bm{C}}_{t}\bm{\eta}_k^s.
\end{equation}
Now the cost term at the instant $t$ in $J_k$ is calculated as
\begin{equation}\label{}
\begin{aligned}
     &\mathbb{E}_{\{\bm{w}_k^s, \bm{v}_k^s,{e}_k^s\}}\{\bm{x}^T_{t+1}\bm{Q}_{t+1}\bm{x}_{t+1} + \bm{u}^T_{t}\bm{R}_{t}\bm{u}_{t}\}\\
     =& \hat{\bm{x}}_{k}^T\tilde{\bm{A}}_{t}^T\bm{Q}_{t+1}\tilde{\bm{A}}_{t}\hat{\bm{x}}_{k} +\bm{U}_{k}^T\tilde{\bm{B}}_{t}^T\bm{Q}_{t+1}\tilde{\bm{B}}_{t}\bm{U}_{k}\\
     &+(\bm{\eta}_k^b)^T\tilde{\bm{C}}_{t}^T\bm{Q}_{t+1}\tilde{\bm{C}}_{t}\bm{\eta}_k^b +2\hat{\bm{x}}_{k}^T\tilde{\bm{A}}_{t}^T\bm{Q}_{t+1}\tilde{\bm{B}}_{t}\bm{U}_{k}\\
     &+2\bm{U}_{k}^T\tilde{\bm{B}}_{t}^T\bm{Q}_{t+1}\tilde{\bm{C}}_{t}\bm{\eta}_k^b +2\hat{\bm{x}}_{k}^T\tilde{\bm{A}}_{t}^T\bm{Q}_{t+1}\tilde{\bm{C}}_{t}\bm{\eta}_k^b\\
     &+\bm{u}^T_{t}\bm{R}_{t}\bm{u}_{t} +\tr\{\tilde{\bm{A}}_{t}^T\bm{Q}_{t+1}\tilde{\bm{A}}_{t}\bm{P}_{k}\}\\
     &+\tr\{\bm{C}_{t}^T\bm{Q}_{t+1}\bm{C}_{t} \diag(\bm{P}_{k}^w, \cdots, \bm{P}_{k+N-1}^w) \}.\\
\end{aligned}
\end{equation}
Then, the expectation of cost $ J_k $ is 
\begin{equation}\label{}
\begin{aligned}
&\mathbb{E}_{\{\bm{w}_k^s, \bm{v}_k^s,{e}_k^s\}} \left\{ J_k  \right\}= \hat{\bm{x}}_{k}^T \left(\sum_{t=k}^{k+N-1}\tilde{\bm{A}}_{t}^T\bm{Q}_{t+1}\tilde{\bm{A}}_{t} \right)\hat{\bm{x}}_{k}\\
&+\bm{U}_{k}^T \left( \sum_{t=k}^{k+N-1} \tilde{\bm{B}}_{t}^T\bm{Q}_{t+1}\tilde{\bm{B}}_{t}\right) \bm{U}_{k} \\
&+ (\bm{\eta}_k^b)^T \left(\sum_{t=k}^{k+N-1} \tilde{\bm{C}}_{t}^T\bm{Q}_{t+1}\tilde{\bm{C}}_{t}\right)\bm{\eta}_k^b \\
&+2\hat{\bm{x}}_{k}^T\left( \sum_{t=k}^{k+N-1} \tilde{\bm{A}}_{t}^T\bm{Q}_{t+1}\tilde{\bm{B}}_{t}\right) \bm{U}_{k} \\
&+2\bm{U}_{k}^T\left( \sum_{t=k}^{k+N-1} \tilde{\bm{B}}_{t}^T\bm{Q}_{t+1}\tilde{\bm{C}}_{t}\right)\bm{\eta}_k^b \\
&+2\hat{\bm{x}}_{k}^T\left( \sum_{t=k}^{k+N-1} \tilde{\bm{A}}_{t}^T\bm{Q}_{t+1}\tilde{\bm{C}}_{t}\right)\bm{\eta}_k^b \\
&+\bm{U}_{k}^T\diag(\bm{R}_k, \cdots, \bm{R}_{k+N-1})\bm{U}_{k} \\
&+\sum_{t=k}^{k+N-1} \tr\{\tilde{\bm{A}}_{t}^T\bm{Q}_{t+1}\tilde{\bm{A}}_{t}\bm{P}_{k}\} \\
&+ \sum_{t=k}^{k+N-1} \tr\{\bm{C}_{t}^T\bm{Q}_{t+1}\bm{C}_{t} \diag(\bm{P}_{k}^w, \cdots, \bm{P}_{k+N-1}^w) \}
\end{aligned}
\end{equation}
Thus the expectation of cost $J_k$ is written in the form stated above, with the following definitions
\begin{equation}\label{Para_A}
\mathscr{A}_k = \sum_{t=k}^{k+N-1}\tilde{\bm{A}}_{t}^T\bm{Q}_{t+1}\tilde{\bm{A}}_{t} ,
\end{equation}
\begin{equation}\label{Para_B}
\mathscr{B}_k = \sum_{t=k}^{k+N-1} \tilde{\bm{B}}_{t}^T\bm{Q}_{t+1}\tilde{\bm{B}}_{t} +\diag(\bm{R}_k, \cdots, \bm{R}_{k+N-1}), 
\end{equation}
\begin{equation}\label{Para_C}
\mathscr{C}_k = \sum_{t=k}^{k+N-1} \tilde{\bm{C}}_{t}^T\bm{Q}_{t+1}\tilde{\bm{C}}_{t}, 
\end{equation}
\begin{equation}\label{Para_D}
\mathscr{D}_k = \sum_{t=k}^{k+N-1} \tilde{\bm{B}}_{t}^T\bm{Q}_{t+1}\tilde{\bm{C}}_{t},
\end{equation}
\begin{equation}\label{Para_b}
\hat{\bm{b}}_k = \left( \sum_{t=k}^{k+N-1} \tilde{\bm{B}}_{t}^T\bm{Q}_{t+1}\tilde{\bm{A}}_{t}\right)\hat{\bm{x}}_{k},
\end{equation}
\begin{equation}\label{Para_c}
\hat{\bm{c}}_k = \left( \sum_{t=k}^{k+N-1} \tilde{\bm{C}}_{t}^T\bm{Q}_{t+1}\tilde{\bm{A}}_{t}\right)\hat{\bm{x}}_{k},
\end{equation}
\begin{equation}
\begin{aligned}
\const &= \sum_{t=k}^{k+N-1} \tr\{\tilde{\bm{A}}_{t}^T\bm{Q}_{t+1}\tilde{\bm{A}}_{t}\bm{P}_{k}\} \\
&+ \sum_{t=k}^{k+N-1} \tr\{\bm{C}_{t}^T\bm{Q}_{t+1}\bm{C}_{t} \diag(\bm{P}_{k}^w, \cdots, \bm{P}_{k+N-1}^w)\}.
\end{aligned}
\end{equation}
\end{proof}

\begin{myTheo}
Problem $(\mathcal{P}_k)$ can be solved by the following semidefinite programming (SDP):
\begin{equation}\label{Prob_k_SDP}
\begin{aligned}
    (\mathcal{P}_k): \quad & \min_{\bm{y}_k} \quad\rho_k \\
    \quad s.t. \quad & \left[ \begin{array}{ccc} \bm{I} & \bm{y}_k & \bm{F}_k \\ \bm{y}_k^T & \rho_k-\tau_1-\tau_2 & -\bm{h}_k^T \\ \bm{F}_k^T &  -\bm{h}_k & \bm{G}_k \end{array} \right]   \geq 0,
\end{aligned}
\end{equation}
in decision variable $\bm{y}_k$, $\rho_k$, $\tau_1$, and $\tau_2$, where
\begin{equation}\label{h_vector}
    \bm{h}_k = \hat{\bm{c}}_k - \mathscr{D}_k^T\mathscr{B}_k^{-1}\hat{\bm{b}}_k,
\end{equation}
\begin{equation}\label{F_matrix}
    \bm{F}_k = \mathscr{B}_k^{-1/2}\mathscr{D}_k.
\end{equation}
\begin{equation}
     \bm{G}_k = -\mathscr{C}_k+ \tau_1\bm{M}_k^{1} + \tau_2\bm{M}_k^2 +\bm{F}_k^T\bm{F}_k,
\end{equation}
\begin{equation}
    \bm{M}_k^{1} =  \left[ \begin{array}{cc} \bm{M}_k^{-1} & 0 \\ 0 & 0 \end{array} \right],
\end{equation}
\begin{equation}
    \bm{M}_k^{2} =  \left[ \begin{array}{cc} 0 & 0 \\ 0 & (\diag(\bm{M}_k^w,\cdots,\bm{M}_{k+N-1}^w))^{-1} \end{array} \right],
\end{equation}
so that the control law is the first component of 
\begin{equation}\label{U_vecotr}
    \bm{U}_k = \mathscr{B}_k^{-1/2}\bm{y}_k-\mathscr{B}_k^{-1}\hat{\bm{b}}_k.
\end{equation}
\end{myTheo}
\begin{proof}
Substituting the control law \eqref{U_vecotr} into cost expectation \eqref{J_K_exp} yields
\begin{equation}\label{J_k_exp_2}
        \mathbb{E}\{J_k\} = \bm{y}_k^T\bm{y}_k+2\bm{h}^T_k\bm{\eta}_k +2\bm{y}_k^T\bm{F}_k\bm{\eta}_k+\bm{\eta}_k^T\mathscr{C}_k\bm{\eta}_k+const.
\end{equation}
With \eqref{J_k_exp_2}, the problem in \eqref{Prob_k} can be reformed as 
\begin{equation}
    \begin{aligned}\label{eq: reformed problem}
        (\mathcal{P}_k):   & \min_{\bm{y}_k} \max_{ \{\bm{e}_k^b, \bm{w}_k^b \}}  \bm{y}_k^T\bm{y}_k+2\bm{h}^T_k\bm{\eta}_k +2\bm{y}_k^T\bm{F}_k\bm{\eta}_k+\bm{\eta}_k^T\mathscr{C}_k\bm{\eta}_k \\
         s.t. \quad &\bm{e}_k^b \in \mathcal{E}(0,\bm{M}_k),\\
         & \bm{W}_k^b \in \mathcal{E}(0,\diag(\bm{M}_k^w,\cdots,\bm{M}_{k+N-1}^w)).
    \end{aligned}
\end{equation}
We introduce the variable $\rho_k$ and rewrite \eqref{eq: reformed problem} as
\begin{equation}\label{prob_roh}
    \begin{aligned}
        (\mathcal{P}_k):  \quad & \min_{\bm{y}_k} \quad \rho_k \\
        \quad s.t. \quad & \rho_k-\bm{y}_k^T\bm{y}_k-2\bm{h}^T_k\bm{\eta}_k -2\bm{y}_k^T\bm{F}_k\bm{\eta}_k -\bm{\eta}_k^T\mathscr{C}_k\bm{\eta}_k \geq 0 ,\\
        \quad & (\bm{e}_k^b)^T \bm{M}_k^{-1} \bm{e}_k^b \leq 1 , \\
        \quad & (\bm{W}_k^b)^T\diag(\bm{M}_k^w,\cdots,\bm{M}_{k+N-1}^w)^{-1}\bm{W}_k^b \leq 1 .
    \end{aligned}
\end{equation}
The first constraint in \eqref{prob_roh} can be reformed as 
\begin{equation}\label{constraint1}
	\begin{aligned}
		\left[ \begin{array}{c} 1 \\ \bm{\eta}_k \end{array} \right]^T \left[ \begin{array}{cc} \rho_k-\bm{y}_k^T\bm{y}_k & -\bm{h}_k^T-\bm{y}_k^T\bm{F}_k \\ -\bm{h}_k-\bm{F}_k^T\bm{y}_k & -\mathscr{C}_k \end{array} \right] \left[ \begin{array}{c} 1 \\ \bm{\eta}_k \end{array} \right] \geq 0.
	\end{aligned}
\end{equation}
The second constraint can be rewritten as
\begin{equation}\label{constraint2}
    \begin{aligned}
        \left[ \begin{array}{c} 1 \\ \bm{\eta}_k \end{array} \right]^T \left[ \begin{array}{cc} 1 & 0 \\ 0 & -\bm{M}_k^{1} \end{array} \right] \left[ \begin{array}{c} 1 \\ \bm{\eta}_k \end{array} \right] \geq 0 \;.
    \end{aligned}
\end{equation}
The third constraint can be reformulated as
\begin{equation}\label{constraint3}
	\begin{aligned}
		\left[ \begin{array}{c} 1 \\ \bm{\eta}_k \end{array} \right]^T \left[ \begin{array}{cc} 1 & 0 \\ 0 & -\bm{M}_k^2 \end{array} \right] \left[ \begin{array}{c} 1 \\ \bm{\eta}_k \end{array} \right] \geq 0.
	\end{aligned}
\end{equation}
According to the S-procedure \cite{boyd1994linear}, for all $\bm{\xi}_t=[1, \bm{\eta}_k]^T$ that satisfies the constraints in (\ref{constraint2}) and (\ref{constraint3}), the constraint in (\ref{constraint1}) also holds if there exist $\tau_1\geq 0$ and $\tau_2\geq 0$ such that 
\begin{equation}\label{LMI_spro}
    \begin{aligned}
        &\left[ \begin{array}{cc} \rho_k-\bm{y}_k^T\bm{y}_k & -\bm{h}_k^T-\bm{y}_k^T\bm{F}_k \\ -\bm{h}_k-\bm{F}_k^T\bm{y}_k & -\mathscr{C}_k \end{array} \right] - \tau_1 \left[ \begin{array}{cc} 1 & 0 \\ 0 & -\bm{M}_k^{1} \end{array} \right] \\
        &- \tau_2 \left[ \begin{array}{cc} 1 & 0 \\ 0 & -\bm{M}_k^2 \end{array} \right]  \geq 0.
    \end{aligned}
\end{equation}
Collecting all terms in (\ref{LMI_spro}) and we get
\begin{equation}\label{LMI_Schur}
	\begin{aligned}
		\left[ \begin{array}{cc} \rho_k-\tau_1-\tau_2 & -\bm{h}_k^T \\ -\bm{h}_k & -\mathscr{C}_k+ \tau_1\bm{M}_k^{1} + \tau_2\bm{M}_k^2 +\bm{F}_k^T\bm{F}_k\end{array} \right] \\
  -  \left[ \begin{array}{cc} \bm{y}_k & \bm{F}_k  \end{array} \right]^T \left[ \begin{array}{cc}  \bm{y}_k & \bm{F}_k \end{array} \right]   \geq 0.
	\end{aligned}
\end{equation}
According to the Schur complement theorem \cite{boyd1994linear}, the inequality can be rewritten as the constraint in \eqref{Prob_k_SDP}, and then the proof is completed.
\end{proof}

\section{Simulations and results} \label{section5}
In this section, we demonstrate the performance of our approach by numerical simulations. We compare our mixed state estimation (MIX) with the Kalman Filter (KF) and ellipsoid set-membership filter (ESM)~\cite{maksarov1996state}, and also compare our proposed robust control with mixed estimation (RCMIX) with the robust control with Kalman filter (RCKF) and ellipsoid set-membership filter (RCESM)~\cite{10591250} independently. 

The simulation considers the control problem $(\mathcal{P})$ in (\ref{Prob1}) with the following particulars
\begin{equation}
	\begin{aligned}
		\bm{A}_k=(1+0.1\sin(k)) \left[ \begin{array}{cc} 0.6 & 0.7 \\ 0.25 & 0.5 \end{array} \right],
	\end{aligned}
\end{equation}
\begin{equation}
	\begin{aligned}
		\bm{B}_k=\left[ \begin{array}{cc} 1  & 0.3  \end{array} \right]^T, \quad
		\bm{C}_k=\left[ \begin{array}{cc} 0.2 & 1  \end{array} \right],
	\end{aligned}
\end{equation}
\begin{equation}
\begin{aligned}
    \bm{w}_k^s \sim \mathcal{N}(0,0.25),\quad
    \bm{w}_k^b \in \mathcal{E}(0,[5 \quad 2; \quad 2 \quad 5]),\\
\end{aligned}
\end{equation}
\begin{equation}
\begin{aligned}
    \bm{v}_k^s \sim \mathcal{N}(0,0.25),\quad
    \bm{v}_k^b \in \mathcal{E}(0,5),\\
\end{aligned}
\end{equation}
\begin{equation}
	\begin{aligned}
		\bm{Q}_k=\left[ \begin{array}{cc} 10 & 0 \\ 0 & 1 \end{array} \right]
, \quad \bm{R}_k=1, \quad N = 100.
	\end{aligned}
\end{equation}
The initial state $\bm{x}_0\sim \mathcal{X}(\hat{\bm{x}}_0, \bm{P}_0, \bm{M}_0)$ is set as
\begin{equation}
	\begin{aligned}
		\hat{\bm{x}}_0= \left[ \begin{array}{c} 60 \\ -45 \end{array} \right],  \bm{P}_0 = \left[ \begin{array}{cc} 10^2 & 0 \\ 0 & 10^2 \end{array} \right], \bm{M}_0 = \left[ \begin{array}{cc} 20^2 & 0 \\ 0 & 20^2 \end{array} \right].
	\end{aligned}
\end{equation}
In this simulation, we assume the bounded process noise $\bm{w}_k^b$ is non-symmetrically distributed, with $90\%$ of normalized noise in each coordinate falling uniformly in $(0,1)$ and $10\%$ in $(0,-1)$. This assumption may arise when the process noise contains systematic components due to neglected dynamics, model parameter errors, or unknown but bounded input.

\begin{figure}[htbp]
	\centering
	\includegraphics[width=0.5\textwidth]{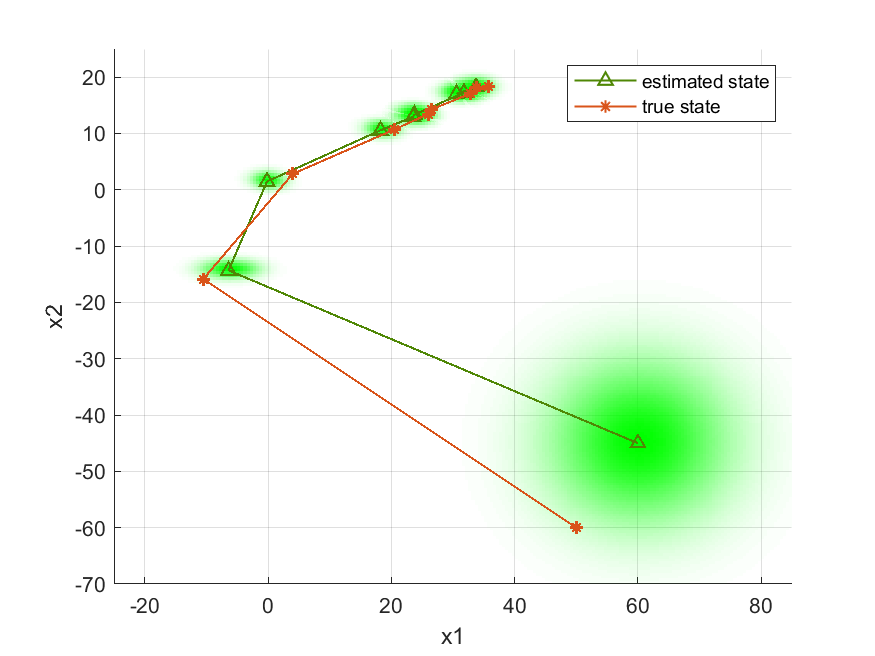}\\
	\caption{State estimation using Kalman filter (KF). } \label{KalmanFilter}
\end{figure}

\begin{figure}[htbp]
	\centering
	\includegraphics[width=0.5\textwidth]{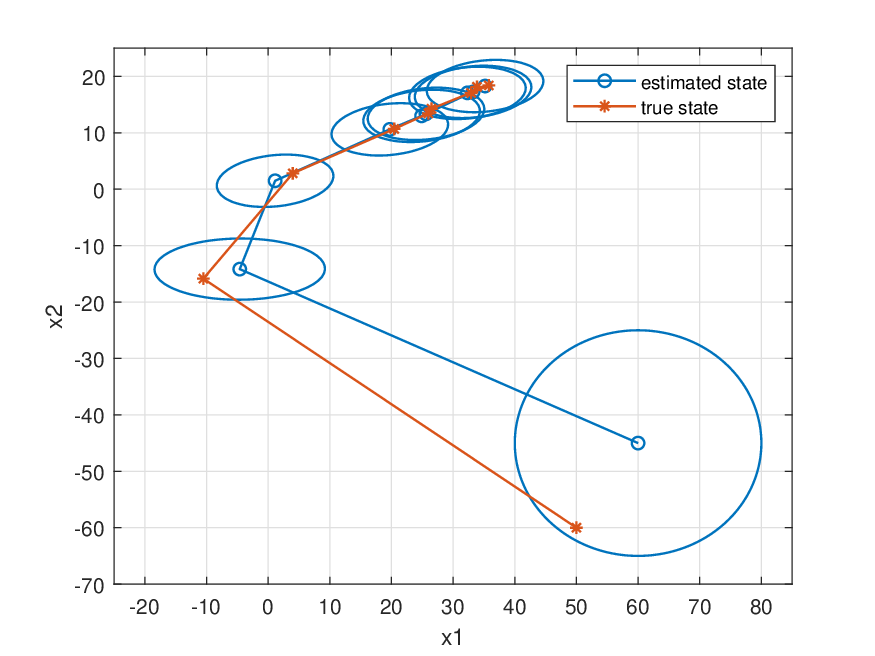}\\
	\caption{State estimation using ellipsoid set-membership filter (ESM). } \label{EllipsoidEst}
\end{figure}

\begin{figure}[htbp]
	\centering
	\includegraphics[width=0.5\textwidth]{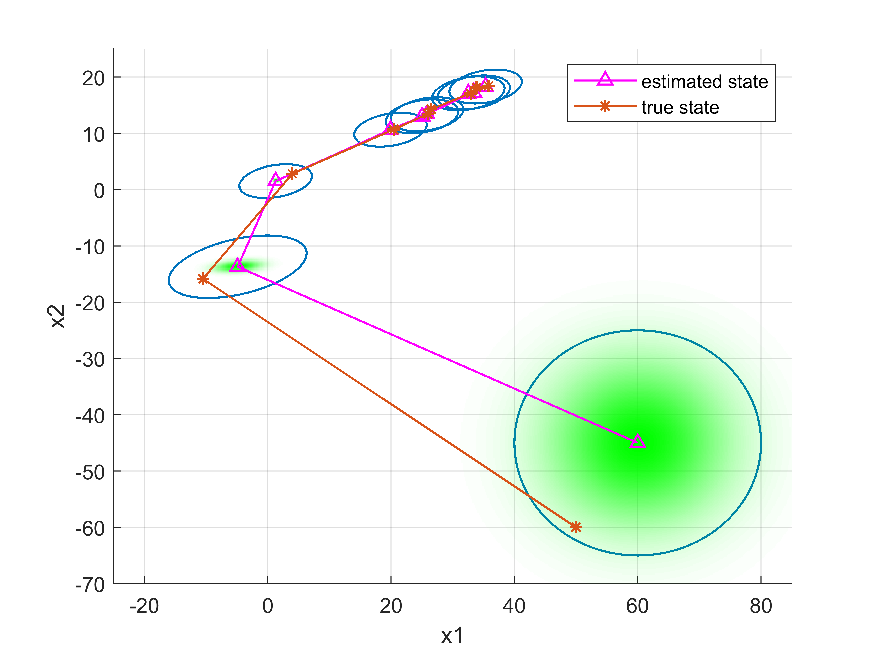}\\
	\caption{State estimation using mixed estimator (MIX). } \label{MixEst}
\end{figure}

Fig.1, Fig.2, and Fig.3 show the state estimation results using KF, ESM, and MIX, respectively. The updated covariance by KF (green area in Fig. 1) is smaller than the updated ellipsoid by MIX (ellipsoid with blue edge in Fig.3), demonstrating that the estimation result by Kalman filter is over-optimistic. The updated ellipsoid by ESM (ellipsoid with blue edge in Fig.2) is larger than the updated ellipsoid by MIX, demonstrating that the estimation result by ellipsoid set-membership filter is over-conservative. The estimated state value by ESM (blue line with circle mark) and MIX (pink line with triangle mark) is closer than KF (green line with triangle mark), reflecting that Kalman filter is sensitive to the non-symmetric process noise, while the ellipsoid set-membership filter and our mixed estimation are less vulnerable.

%使用均方误差（MSE）、均方根误差（RMSE）、平均绝对误差（MAE）% Confidence level = 95%
\begin{table}[H]
	\caption{Comparison between estimation performance with Kalman filter(KF), ellipsoid set-membership filter(ESM), mixed estimator(MIX).}
	\centering
	\resizebox{6cm}{!}
	{\begin{tabular}{ccccc}
			\toprule
			Estimator &  KF  & ESM & MIX  \\ 
			\midrule
			MAE       & 1.5689 & 1.0636 & 1.0478 \\ 
			MSE       & 4.8169 & 3.0427 & 2.9832 \\
			RMSE      & 2.1917 & 1.7432 & 1.7262 \\
			Volume    & 121.2270 & 1869.2 & 293.0740 \\
			Trace     & 26.0975 & 111.2582 & 44.1025\\
			\bottomrule
	\end{tabular}}
	\label{table1}
\end{table}

TABLE \ref{table1} shows the mean absolute error (MAE), mean squared error (MSE), root mean squared error (RMSE), the volume of the estimated ellipsoid set (the ellipsoid with the $95\%$ confidence level of covariance estimated by KF), the trace of ellipsoid shape matrix for the state estimation results by KF, ESM and MIX, respectively. The MAE, MSE and RMSE for our mixed estimation are the smallest among the three approaches, demonstrating that our method has the least estimation error. The value for ellipsoid shape matrix volume and trace for our mixed estimation is greater than those of KF and smaller than those of ESM, meaning that the estimation result of MIX is not over-optimistic as KF and not over-conservative as ESM.

\begin{figure}[htbp]
	\centering
	\includegraphics[width=0.5\textwidth]{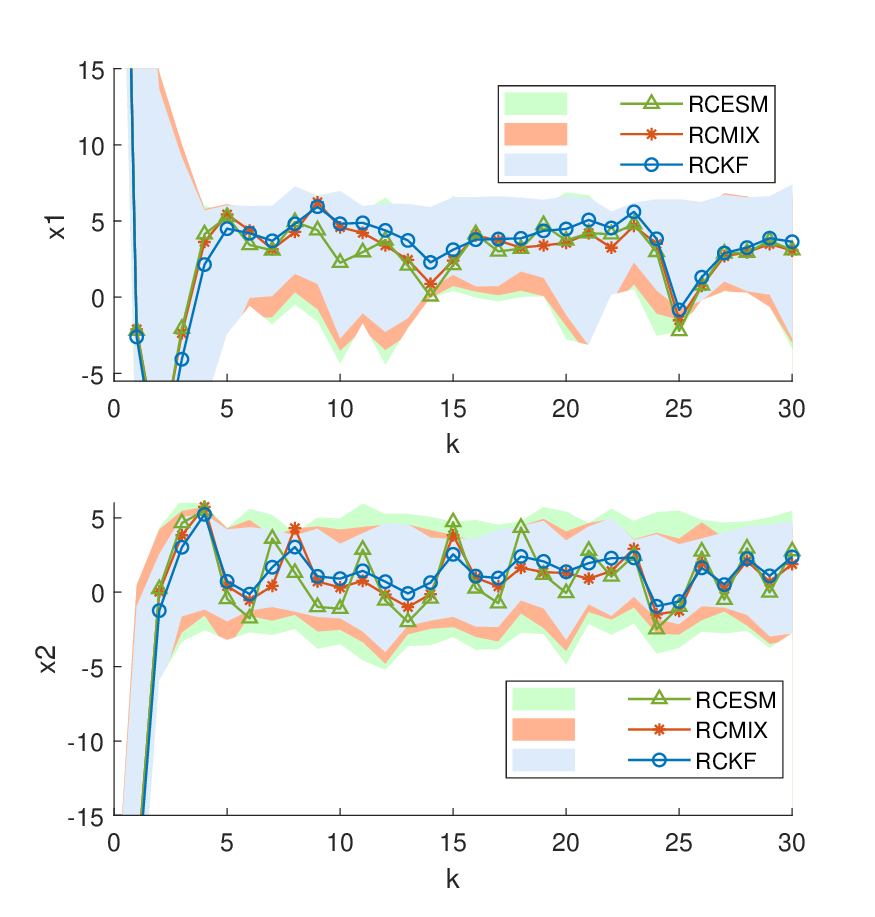}\\
	\caption{State regulating under RCKF, RCESM and RCMIX. } \label{stateregualtion}
\end{figure}

\begin{table}[H]
\caption{Comparison between control performance with robust control with Kalman filter(RCKF), robust control with ellipsoid set-membership filter(RCESM), and robust control with mixed estimator(RCMIX).}
\centering
\resizebox{6cm}{!}
{\begin{tabular}{ccccc}
        \toprule
        Controller &  RCKF  &  RCESM &  RCMIX  \\ 
        \midrule
        MAE  &  3.3879 &  3.0334 &  2.9240 \\ 
        MSE  & 46.3215 & 43.4079 & 42.9013 \\
        RMSE &  6.8060 &  6.5885 &  6.5499 \\
        \bottomrule
\end{tabular}}
\label{table2}
\end{table}

Fig. 4 shows the 100 rounds of simulation results of state regulating control with RCKF, RCESM and RCMIX, respectively. The RCESM has the widest green area among the three, demonstrating that control results of RCESM is over-conservative. The RCKF has the most narrowed blue area with the least conservativeness, while having a larger deviation from the reference value zero than RCMIX. TABLE \ref{table2} shows that our RCMIX has the smallest state regulating error among three methods.

\section{Conclusion} \label{section6}
We present a linear quadratic control approach for systems where the process and observation noises are unknown, and the system states cannot be measured directly. To account for both types of noise—stochastic and bounded—we incorporate a Gaussian distribution and an ellipsoidal set for modeling system noise in the control design. Our approach integrates Kalman and set-membership filters to design a state estimation method. The designed estimator is used to approximate the unmeasurable system states and to update the covariance and shape matrix of the stochastic and bounded estimation errors. This results in estimates that are less optimistic than those of the Kalman filter while being less conservative than those of the set-membership filter. Based on this state estimation, we derive a robust control law by jointly minimizing the performance metric for stochastic control and the worst-case control metric. The resulting control law mitigates over-conservativeness while ensuring stability in the presence of both stochastic and bounded noises. Future research will focus on reducing the computational complexity of the control design, as the semidefinite programming involved in deriving the control law for noisy systems can be computationally expensive, particularly for large system states.

\bibliographystyle{unsrt}
\bibliography{bibfile}

\begin{thebibliography}{10}

\bibitem{11219403}
Yajie Yu, Xuehui Ma, Shiliang Zhang, Zhuzhu Wang, Xubing Shi, Yushuai Li, and Tingwen Huang.
\newblock Adaptive ensemble control for stochastic systems with mixed asymmetric laplace noises.
\newblock {\em IEEE Transactions on Systems, Man, and Cybernetics: Systems}, pages 1--14, 2025.

\bibitem{ma2022adaptivequantile}
Xuehui Ma, Fucai Qian, Shiliang Zhang, and Li~Wu.
\newblock Adaptive quantile control for stochastic system.
\newblock {\em ISA transactions}, 123:110--121, 2022.

\bibitem{10591250}
Xuehui Ma, Yutong Chen, Shiliang Zhang, Yushuai Li, Fucai Qian, and Zhiyong Sun.
\newblock Robust quadratic optimal control of linear systems with ellipsoid-set learning.
\newblock In {\em 2024 European Control Conference (ECC)}, pages 2125--2131, 2024.

\bibitem{lorenzen2019robust}
Matthias Lorenzen, Mark Cannon, and Frank Allgöwer.
\newblock Robust mpc with recursive model update.
\newblock {\em Automatica}, 103:461--471, 2019.

\bibitem{aastrom2012introduction}
Karl~J {\AA}str{\"o}m.
\newblock {\em Introduction to Stochastic Control Theory}.
\newblock Courier Corporation, 2012.

\bibitem{bertsimas2007constrained}
Dimitris Bertsimas and David~B Brown.
\newblock Constrained stochastic {LQC}: a tractable approach.
\newblock {\em IEEE Transactions on Automatic Control}, 52(10):1826--1841, 2007.

\bibitem{di2004set}
Mauro Di~Marco, Andrea Garulli, Antonio Giannitrapani, and Antonio Vicino.
\newblock A set theoretic approach to dynamic robot localization and mapping.
\newblock {\em Autonomous robots}, 16(1):23--47, 2004.

\bibitem{hanebeck1999new}
Uwe~D Hanebeck and Joachim Horn.
\newblock A new estimator for mixed stochastic and set theoretic uncertainty models applied to mobile robot localization.
\newblock In {\em Proceedings 1999 IEEE International Conference on Robotics and Automation}, volume~2, pages 1335--1340. IEEE, 1999.

\bibitem{liu2025adaptive}
Haolin Liu, Shiliang Zhang, Xiaohui Zhang, Shangbin Jiao, Xuehui Ma, Ting Shang, Yan Yan, Wenqi Bai, and Youmin Zhang.
\newblock Adaptive lyapunov-constrained mpc for fault-tolerant auv trajectory tracking.
\newblock {\em arXiv preprint arXiv:2509.17237}, 2025.

\bibitem{SOLIMAN2022387}
Mohamed Soliman, Bruno Morabito, and Rolf Findeisen.
\newblock Towards safe exploration for autonomous vehicles using dual model predictive control.
\newblock In {\em 9th IFAC Symposium on Mechatronic Systems MECHATRONICS 2022}, volume~55, pages 387--392, 2022.

\bibitem{yang2022set}
Hao Yang, Yilian Zhang, Wei Gu, Fuwen Yang, and Zhiquan Liu.
\newblock Set-membership filtering for automatic guided vehicles with unknown-but-bounded noises.
\newblock {\em Transactions of the Institute of Measurement and Control}, 44(3):716--725, 2022.

\bibitem{ma2025fault}
Xuehui Ma, Shiliang Zhang, and Zhiyong Sun.
\newblock Fault tolerant control of mecanum wheeled mobile robots.
\newblock {\em arXiv preprint arXiv:2512.06444}, 2025.

\bibitem{liu2025adaptiveFault}
Haolin Liu, Shiliang Zhang, Shangbin Jiao, Xiaohui Zhang, Xuehui Ma, Yan Yan, Wenchuan Cui, and Youmin Zhang.
\newblock Adaptive fault-tolerant control of underwater vehicles with thruster failures.
\newblock {\em arXiv preprint arXiv:2504.16037}, 2025.

\bibitem{MISHRA2021109512}
Prabhat~K. Mishra, Sanket~S. Diwale, Colin~N. Jones, and Debasish Chatterjee.
\newblock Reference tracking stochastic model predictive control over unreliable channels and bounded control actions.
\newblock {\em Automatica}, 127:109--512, 2021.

\bibitem{arcari2023stochastic}
Elena Arcari, Andrea Iannelli, Andrea Carron, and Melanie~N Zeilinger.
\newblock Stochastic mpc with robustness to bounded parameteric uncertainty.
\newblock {\em IEEE Transactions on Automatic Control}, 68(12):7601--7615, 2023.

\bibitem{ma2022active}
Xuehui Ma, Shiliang Zhang, Fucai Qian, Jinbao Wang, and Yushuai Li.
\newblock Active learning for anti-disturbance dual control of unknown nonlinear systems.
\newblock {\em arXiv preprint arXiv:2212.08934}, 2022.

\bibitem{de2024state}
Alesi~A de~Paula, Guilherme~V Raffo, and Bruno~OS Teixeira.
\newblock State estimators for discrete-time descriptor linear systems with mixed uncertainties and state constraints.
\newblock {\em International Journal of Robust and Nonlinear Control}, 34(12):7936--7967, 2024.

\bibitem{lu2017multi}
Kelin Lu and Rui Zhou.
\newblock Multi-sensor fusion for robust target tracking in the simultaneous presence of set-membership and stochastic gaussian uncertainties.
\newblock {\em IET Radar, Sonar \& Navigation}, 11(4):621--628, 2017.

\bibitem{iannelli2020structured}
Andrea Iannelli, Mohammad Khosravi, and Roy~S Smith.
\newblock Structured exploration in the finite horizon linear quadratic dual control problem.
\newblock {\em IFAC-PapersOnLine}, 53(2):959--964, 2020.

\bibitem{parsi2022explicit}
Anilkumar Parsi, Andrea Iannelli, and Roy~S Smith.
\newblock An explicit dual control approach for constrained reference tracking of uncertain linear systems.
\newblock {\em IEEE Transactions on Automatic Control}, 68(5):2652--2666, 2022.

\bibitem{parsi2022scalable}
Anilkumar Parsi, Andrea Iannelli, and Roy~S Smith.
\newblock Scalable tube model predictive control of uncertain linear systems using ellipsoidal sets.
\newblock {\em International Journal of Robust and Nonlinear Control}, pages 1--22, 2022.

\bibitem{10551447}
Xuehui Ma, Shiliang Zhang, Yushuai Li, Fucai Qian, Zhiyong Sun, and Tingwen Huang.
\newblock Adaptive robust tracking control with active learning for linear systems with ellipsoidal bounded uncertainties.
\newblock {\em IEEE Transactions on Automatic Control}, 69(11):8096--8103, 2024.

\bibitem{10904011}
Marco Casini, Andrea Garulli, and Antonio Vicino.
\newblock Adaptive threshold selection for set membership state estimation with quantized measurements.
\newblock {\em IEEE Transactions on Automatic Control}, pages 1--14, 2025.

\bibitem{noack2012combined}
Benjamin Noack, Florian Pfaff, and Uwe~D Hanebeck.
\newblock Combined stochastic and set-membership information filtering in multisensor systems.
\newblock In {\em 2012 15th International Conference on Information Fusion}, pages 1218--1224. IEEE, 2012.

\bibitem{maksarov1996state}
DG~Maksarov and JP~Norton.
\newblock State bounding with ellipsoidal set description of the uncertainty.
\newblock {\em International Journal of Control}, 65(5):847--866, 1996.

\bibitem{noack2012optimal}
Benjamin Noack, Florian Pfaff, and Uwe~D Hanebeck.
\newblock Optimal kalman gains for combined stochastic and set-membership state estimation.
\newblock In {\em 2012 IEEE 51st IEEE Conference on Decision and Control (CDC)}, pages 4035--4040. IEEE, 2012.

\bibitem{boyd1994linear}
Stephen Boyd, Laurent El~Ghaoui, Eric Feron, and Venkataramanan Balakrishnan.
\newblock {\em Linear Matrix Inequalities in System and Control Theory}.
\newblock SIAM, 1994.

\end{thebibliography}

\end{document}